\documentclass{article}
\usepackage[utf8]{inputenc}

\usepackage{amsmath}
\usepackage{amssymb}
\usepackage{amsthm}

\AtEndDocument{\noindent DEPARTMENT OF MATHEMATICAL SCIENCES, SEOUL NATIONAL UNIVERSITY, SEOUL, 151-747, KOREA. E-mail adress: sagelee24@naver.com
}

\newtheorem{theorem}{Theorem}[section]
\newtheorem{lemma}[theorem]{Lemma}
\newtheorem{corollary}[theorem]{Corollary}

\newtheorem{remark}[theorem]{Remark}
\newtheorem{definition}[theorem]{Definition}
\newtheorem{assumption}[theorem]{Assumption}
\newtheorem{notation}[theorem]{Notation}

\numberwithin{equation}{section}

\newcommand{\vertiii}[1]{{\left\vert\kern-0.25ex\left\vert\kern-0.25ex\left\vert #1 
    \right\vert\kern-0.25ex\right\vert\kern-0.25ex\right\vert}}

\title{The Invariant Subspace Problem, Revision of arXiv:2002.11533v11 [math.GM]}
\author{SA GE LEE}

\date{December 2022}

\begin{document}

\maketitle

\begin{abstract}
The invariant subspace problem is solved correcting my earlier attempts [6]--[12].
\end{abstract}

\section{Introduction}\label{sec:1}
Let $T$ be an arbitrarily chosen bounded linear operator on a separable infinite dimensional Hilbert space $H$ over $\mathbb{C}$ or $\mathbb{R}$.

The invariant subspace problem: Does $T$ have a nontrivial closed invariant subspace of $H$ for $T$? ([6], pp.95--96, \S 151).

As far as the invariant subspace problem is concerned, we may and shall assume:

\begin{assumption}\label{assumption:1.1}
$T$ is contractive ($\|T\| \leq 1$) and admits a unit cyclic vector $e_1$ ($\|e_1\| =1$, $ H = \left[T^{n-1}e_1 : n \in \mathbb{N} \right]$, where $[ \cdots ] $ denotes the closed linear span).
\end{assumption}

\begin{notation}\label{notation:1.2}
By the Gram-Schmidt orthonormalization, we obtatin the orthonormal basis of $H$:

\begin{equation}\label{eq:1.1}
e = (e_n : n \in \mathbb{N}),
\end{equation}
satisfying
\begin{equation}\label{eq:1.2}
[e_1, e_2 , \cdots, e_n ] = [e_1 , Te_1 , T^2 e_1 , \cdots , T^{n-1} e_1 ], \quad \forall n \in \mathbb{N}.
\end{equation}

\begin{align}\label{eq:1.3}
B(H) = \textup{The von Neumann algebra of bounded linear operators on $H$, equipped} \\
\textup{with the strong operator topology (SO) and the operator norm $\|\cdot\|$.} \nonumber
\end{align}

Let $\langle \cdot , \cdot \rangle$ be the inner product for $H$ and $\| \cdot \|_e$ be the new norm on $B(H)$ defined by:
\begin{equation}\label{eq:1.4}
\|A\|_e = \sum_{k, l \in \mathbb{N}} \frac{1}{2^{k+l}}  | \langle e_k , A e_l \rangle | ,\quad \forall A \in B(H).
\end{equation}

\begin{equation}\label{eq:1.5}
E_n = \textup{The projection $\in B(H)$, where $\text{Range}(E_n) = [e_1, \cdots, e_n]$, ${}^\forall n \in \mathbb{N}. $}
\end{equation}

\begin{equation}\label{eq:1.6}
\mathcal{A} = \textup{The abelian von Neumann algebra generated by $(E_n : n \in \mathbb{N})$.}
\end{equation}

\begin{equation}\label{eq:1.7}
\mathcal{A}_1 = \{ A \in \mathcal{A} : \| A \| \leq 1 \}
\end{equation}

\begin{align}\label{eq:1.8}
\mathcal{A}_{1,0} = \Big\{ \alpha_1 (E_2 - E_1 ) + \alpha_2 (E_3 - E_2 ) + \alpha_3 (E_4 - E_3) + \cdots + \alpha_i (E_{i+1} - E_i ) + \cdots \\
: \alpha_i \in \mathbb{R}, |\alpha| \leq 1 , \forall i \in \mathbb{N} \Big\}. \nonumber
\end{align}

Note that $\mathcal{A}_{1,0}$ is a WO-closed subset of $\mathcal{A}_1$, consisting of self-adjoint elements of $\mathcal{A}_1$ so that for every $A \in \mathcal{A}_{1,0}$,

\begin{equation}\label{eq:1.9}
(\textup{Kernel} A )^\perp = \overline{ \textup{Range} A}
\end{equation}

Recall that the Banach space $l^\infty = \Big\{ (\alpha_i : \alpha_i \in \mathbb{R}, {}^\forall i \in \mathbb{N} ) \text{ with } \|(\alpha_i : i \in \mathbb{N} ) \|_\infty = \sup_{i \in \mathbb{N}} | \alpha_i | < \infty \Big\}$ is the dual space of the Banach space $l^1 = \Big\{ (\beta_i : \beta_i \in \mathbb{R}, {}^\forall i \in \mathbb{N}) \text{ with } \| (\beta_i: i \in \mathbb{N} ) \|_1 = \sum_{i \in \mathbb{N}} |\beta_i| < \infty \Big\}$.

On the other hand, there exists one to one correspondence between the unit ball $(l^\infty)_1$ of $l^\infty$ and $\mathcal{A}_{1,0}$, via

\begin{equation}\label{eq:1.10}
( \alpha_i : i \in \mathbb{N}) \in (l^\infty)_1 \longleftrightarrow \sum_{i \in \mathbb{N}} \alpha_i (E_{i+1} - E_i ) \in \mathcal{A}_{1, 0}.
\end{equation}

Therefore, we can transplant the weak-$*$ topology $\sigma( l^\infty, l^1 )$ relativised on $(l^\infty)_1$ onto $\mathcal{A}_{1,0}$. Hence $\big(\mathcal{A}_{1,0} , \sigma(l^\infty, l^1)\big)$ is regarded as a separable compact Hausdorff space ([4] p.434 Theorem 2 and p.426 Theorem 1).

For every $n \in \mathbb{N}$, we define:

\begin{align}\label{eq:1.11}
\mathcal{A}(n) =& \Big\{ A \in \mathcal{A}_{1,0} : \langle Ae_j, e_j \rangle = 0, 1 \leq {}^\forall j \leq n, \\
&| \sum_{i \in \mathbb{N}} \langle A e_i, e_i \rangle \beta_i | \geq | \sum_{i \in \mathbb{N}} \langle B_n e_i , e_i \rangle \beta_i |, {}^\forall (\beta_i : i \in \mathbb{N}) \in (l^1 )_1, \nonumber \\
& \text{ such that } \beta_j = 0 , 1 \leq {}^\forall j \leq n \Big\}, \nonumber
\end{align}
where $(l^1)_1$ is the unit ball of $l^1$.

\begin{align}\label{eq:1.12}
B_n = (E_{n+1} - E_n ) + (E_{n+2} - E_{n+1} ) + (E_{n+3} - E_{n+2}) + \\
\cdots + (E_{n+k} - E_{n+k-1} )+ \cdots \quad (k \in \mathbb{N}). \nonumber
\end{align}

\end{notation}

\begin{lemma}\label{lemma:1.3}
${}^\forall n \in \mathbb{N}$,
\begin{equation}\label{eq:1.13}
B_n \in \mathcal{A}(n) , \text{ $\big( \mathcal{A}(n), \sigma (l^\infty , l^1 ) \big)$ is a separable compact Hausdorff space,}
\end{equation}
\begin{equation}\label{eq:1.14}
0 \notin \mathcal{A}(n), \text{ and }
\end{equation}
\begin{equation}\label{eq:1.15}
\mathcal{A}(n) \supset \mathcal{A}(n+1).
\end{equation}
\end{lemma}
\begin{proof}
We omit the easy verification of \eqref{eq:1.13}. (cf. [4] p.434 Theorem 3)

To verify \eqref{eq:1.14}, we observe in \eqref{eq:1.11};
\begin{align*}
\sup \Big\{ | \sum_{i \in \mathbb{N}} \langle A e_i , e_i \rangle \beta_i | : (\beta_i : i \in \mathbb{N}) \in (l^1)_1 \text{ such that } \beta_j = 0, 1 \leq {}^\forall j \leq n \Big \} \hspace{0.5cm} \\
\geq \sup \Big\{ | \sum_{i \in \mathbb{N}} \langle B_n e_i , e_i \rangle \beta_i | : (\beta_i : i \in \mathbb{N}) \in (l^1)_1 \text{ such that } \beta_j = 0, 1 \leq {}^\forall j \leq n \Big \} \\
= \sup \Big\{ | \sum_{i \in \mathbb{N}} \langle B_n e_i , e_i \rangle \beta_i | : (\beta_i : i \in \mathbb{N}) \in (l^1)_1 \Big \} \\
\text{(, since $B_n e_j = 0, \hspace{0.1cm} 1 \leq {}^\forall j \leq n )$} \\
= \|B_n \|_\infty \text{ (where $B_n$ is regarded as an element of $l^\infty$) } = \| B_n \| = 1.
\end{align*}

This gives rise to \eqref{eq:1.14}.

Finally to verify \eqref{eq:1.15}, let
\begin{equation}\label{eq:1.16}
A \in \mathcal{A}(n+1).
\end{equation}

By \eqref{eq:1.11}, where $n$ there, now replaced by $n+1$, we obtain:

\begin{align}\label{eq:1.17}
&A \in \mathcal{A}_{1,0} , \quad \langle A e_j, e_j \rangle = 0, \quad 1 \leq {}^\forall j \leq n+1 \nonumber \\
&\text{ and } \nonumber \\
&\forall (\beta_i: i \in \mathbb{N}) \in (l^1)_1 \text{ such that } \\
&\beta_j = 0, \quad 1 \leq {}^\forall j \leq n+1, \nonumber \\
&| \sum_{i \in \mathbb{N}} \langle A e_i, e_i \rangle \beta_i | \geq | \sum_{i \in \mathbb{N}} \langle B_{n+1} e_i , e_i \rangle \beta_i |. \nonumber
\end{align}

In \eqref{eq:1.11} we notice that
\begin{equation}\label{eq:1.18}
\langle B_n e_i , e_i \rangle = \begin{cases} 1, \quad \text{if } i \geq n+1 \\
0, \quad \text{if } 1 \leq i \leq n \end{cases}
\end{equation}
, in particular, in \eqref{eq:1.17},
\begin{equation}\label{eq:1.19}
\langle B_{n+1} e_i, e_i \rangle = \begin{cases} 1, \quad \text{if } i \geq n+2 \\
0, \quad \text{if } 1 \leq i \leq n+1 \end{cases}
\end{equation}

Thus \eqref{eq:1.17} can be rewritten as follows.
\begin{align}\label{eq:1.20}
&A \in \mathcal{A}_{1,0} , \quad \langle Ae_j, e_j \rangle = 0, \quad 1 \leq {}^\forall j \leq n+1 \nonumber \\
&\text{and} \nonumber \\
&\forall (\beta_i : i \in \mathbb{N}) \in (l^1)_1 \text{ such that } \beta_j = 0, \\
& 1 \leq {}^\forall j \leq n+1, \nonumber \\
&| \sum_{i \in \mathbb{N}} \langle A e_i, e_i \rangle \beta_i | \geq | \sum_{n+2 \leq i < \infty} \beta_i |. \nonumber
\end{align}
, equivalently,
\begin{align}\label{eq:1.21}
&A \in \mathcal{A}_{1,0} , \quad \langle Ae_j, e_j \rangle = 0, \quad 1 \leq {}^\forall j \leq n+1 \nonumber \\
&\text{and} \nonumber \\
&\forall (\beta_i : i \in \mathbb{N}) \in (l^1)_1 \text{ such that } \beta_j = 0, \\
& 1 \leq {}^\forall j \leq n+1, \nonumber \\
&| \sum_{i \in \mathbb{N}} \langle A e_i, e_i \rangle \beta_i | \geq | \sum_{n \leq i < \infty} \beta_i |. \nonumber
\end{align}
, since $\beta_n = \beta_{n+1} = 0$ (cf. \eqref{eq:1.20}).

Now, $\forall(\beta_i: i \in \mathbb{N}) \in (l^1)_1$ such that
\begin{equation}\label{eq:1.22}
\beta_1 = \beta_2 = \cdots = \beta_n = 0,
\end{equation}
we have:
\begin{align}\label{eq:1.23}
|\sum_{i \in \mathbb{N}} \langle B_n e_i , e_i \rangle \beta_i | &= |\sum_{ n+1 \leq i < \infty} \langle B_n e_i , e_i \rangle \beta_i | \\
&= |\sum_{n+1 \leq i < \infty} \beta_i | \quad \text{(, by \eqref{eq:1.18})} \nonumber \\
&= | \sum_{n \leq i <\infty} \beta_i | \quad \text{(, since $\beta_n = 0$ in \eqref{eq:1.22})} \nonumber
\end{align}

Then \eqref{eq:1.21} implies the following by aid of \eqref{eq:1.23}:

\begin{align}\label{eq:1.24}
&A \in \mathcal{A}_{1,0} , \quad \langle Ae_j, e_j \rangle = 0, \quad 1 \leq {}^\forall j \leq n \nonumber \\
&\text{and} \nonumber \\
&\forall (\beta_i : i \in \mathbb{N}) \in (l^1)_1 \text{ such that } \beta_j = 0, \quad 1 \leq {}^\forall j \leq n \\
&\text{we have:} \nonumber \\
&| \sum_{i \in \mathbb{N}} \langle A e_i, e_i \rangle \beta_i | \geq | \sum_{i \in \mathbb{N}} \langle B_n e_i , e_i \rangle \beta_i |. \nonumber
\end{align}
, since $\langle A e_{n+1} , e_{n+1} \rangle = 0$ (cf. \eqref{eq:1.16}), recalling \eqref{eq:1.22} and \eqref{eq:1.23}.

By \eqref{eq:1.24} and \eqref{eq:1.11}, we now can say that our $A \in \mathcal{A}(n+1)$ in \eqref{eq:1.16} also satisfies $A \in \mathcal{A}(n)$. This proves \eqref{eq:1.15}: $\mathcal{A}(n) \supset \mathcal{A}(n+1)$ as desired.
\end{proof}

\section{The solution}\label{sec:2}

\begin{theorem}\label{theorem:2.1}
Let $T$ satisfy Assumption~\ref{assumption:1.1}, without loss of generality. Then $T$ has a nontrivial closed invariant subspace whose relevant projection belongs to $\mathcal{A}$.
\end{theorem}
\begin{proof}
By Lemma~\ref{lemma:1.3}, the family $\{ \mathcal{A}(n) : n \in \mathbb{N} \}$ consists of nonempty closed subsets of the compact Hausdorff space $( \mathcal{A}_{1,0}, \sigma(l^\infty, l^1))$. Hence, by a well known fact ([3] p.223 Theorem 1.3(2)), we obtain that $\bigcap_{n \in \mathbb{N}} \mathcal{A}(n) \neq \emptyset$, so there exists, say,
\begin{equation}\label{eq:2.1}
M \in \bigcap_{n \in \mathbb{N}} \mathcal{A}(n).
\end{equation}

By \eqref{eq:2.1}, \eqref{eq:1.11} and \eqref{eq:1.8}, $M$ can be expressed in many different ways as follows, depending on various ways of choosing $n$ such that $M \in \mathcal{A}(n)$:
\begin{align}\label{eq:2.2}
M = \alpha_{n,n} (E_{n+1}-E_n) + \alpha_{n, n+1} (E_{n+2} - E_{n+1}) + \alpha_{n, n+2} (E_{n+3} - E_{n+2}) + \cdots, \\
\text{where $\alpha_{n,k} \in \mathbb{R}$ such that $| \alpha_{n,k} | \leq 1, \hspace{0.1cm} n \leq {}^\forall k < \infty$.} \nonumber
\end{align}

Now, by \eqref{eq:1.4},
\begin{align}\label{eq:2.3}
\| TM - MTM \|_e = \| (I-M) TM\|_e = \sum_{k,l \in \mathbb{N}} \frac{1}{2^{k+l}} | \langle e_k, (I-M) TM e_l \rangle | \\
\leq \sum_{k, l \in \mathbb{N}} \frac{1}{2^{k+l}} \| (I-M) TM e_l \| \leq \sum_{ l \in \mathbb{N}} \frac{1}{2^{l}} ( 1+ \|M\|) \|T\| \|Me_l\| \nonumber \\
\leq 2 \sum_{l \in \mathbb{N}} \frac{1}{2^l} \|Me_l\| \text{(since $\|M\| \leq 1, \|T\| \leq 1$)} \nonumber \\
= 2 \sum_{l \in \mathbb{N}} \frac{1}{2^l} \| \alpha_{n,n} (E_{n+1} - E_n ) e_l + \alpha_{n, n+1} (E_{n+2} - E_{n+1} ) e_l \nonumber \\ + \alpha_{n, n+2} (E_{n+3} - E_{n+2}) e_l + \cdots \| \nonumber \\
= 2 \sum_{l \in \mathbb{N}} \frac{1}{2^l} \| \alpha_{n,n} \delta_{n+1, l} e_l + \alpha_{n, n+1} \delta_{n+2,l} e_l + \alpha_{n, n+2} \delta_{n+3, l} e_l + \cdots \| \nonumber \\
\text{(where $\delta_{k,l}$ ($k \geq n+1$) are the Kronecker deltas)} \nonumber \\
\leq 2 \sum_{l \in \mathbb{N}} \frac{1}{2^l} \Big( |\alpha_{n,n}| \delta_{n+1, l} + |\alpha_{n, n+1} | \delta_{n+2, l} + | \alpha_{n, n+2} | \delta_{n+3, l } + \cdots \Big) \nonumber \\
\leq 2 \Big( \frac{1}{2^{n+1}} + \frac{1}{2^{n+2}} + \frac{1}{2^{n+3}} + \cdots \Big) \nonumber \\
\text{(, noticing that $|\alpha_{n,k}| \leq 1 , \quad n \leq {}^\forall k < \infty$ in \eqref{eq:2.2}} \nonumber \\
= 2 \cdot \frac{1}{2^n} \rightarrow 0 \quad (n \rightarrow \infty) \nonumber
\end{align}

Hence
\begin{equation}\label{eq:2.4}
\|TM - MTM\|_e = 0
\end{equation}
, i.e.,
\begin{equation}\label{eq:2.5}
TM = MTM.
\end{equation}

By \eqref{eq:2.2},
\begin{equation}\label{eq:2.6}
ME_1 = 0.
\end{equation}

By \eqref{eq:1.14}, \eqref{eq:2.1}, and \eqref{eq:2.2},
\begin{equation}\label{eq:2.7}
M \neq 0.
\end{equation}

by \eqref{eq:2.6} and \eqref{eq:2.7},
\begin{equation*}
\{0 \} \subsetneq \textup{Kernel}(M) \subsetneq H.
\end{equation*}

Hence,
\begin{equation}\label{eq:2.8}
\{ 0 \} \subsetneq \overline{\textup{Range}(M)} \subsetneq H,
\end{equation}

since $\{\textup{Kernel}(M)\}^\perp = \overline{\textup{Range}(M)}$ by \eqref{eq:1.9}. Then, by \eqref{eq:2.8} and \eqref{eq:2.5}, $\overline{\textup{Range}(M)}$ is a nontrivial closed invariant subspace of $H$ for $T$, while the range projection of $M$ belongs to $\mathcal{A}$.
\end{proof}

\begin{remark}\label{remark:2.3}
\begin{enumerate}
\item For the operator $M$ in \eqref{eq:2.1}, we actually have:
\begin{equation}\label{eq:2.9}
\bigcap_{n \in \mathbb{N}} \mathcal{A}(n) = \{ M\},\textup{ a singleton set.} 
\end{equation}
These can be verified by a similar arguments in \eqref{eq:2.4}--\eqref{eq:2.5}. 
\begin{equation}\label{eq:2.10}
{\rm We\ put}\ E =\ {\rm the\ projection\ whose\ range\ is\ }\overline{{\rm Range}\,M}.
\end{equation}

\item Let $F$ be the invariant subspace projection associated with another unit cyclic vector $f_1$ for $T$. Then $UEU^* = F$ for some unitary operator $U$ ([8] Theorem 2.7).
\end{enumerate}
\end{remark}

\begin{definition}\label{definition:2.4}
Let $H$ be a vector space over the quaternion skew-field $\mathbb{H}$, where $\mathbb{H}$ acts on $H$ from the right hand sides of the vectors of $H$. We equip $H$ with the inner product $\langle \cdot, \cdot \rangle: H \times H \rightarrow \mathbb{H}$ of the following sense: $\forall x, y, z \in H$, $\forall \alpha \in \mathbb{H}$,
\begin{equation}\label{eq:2.11}
\begin{cases}
\langle x , y \rangle + \langle x , z \rangle = \langle x , y+z \rangle, \\
\langle x, y \alpha \rangle = \langle x, y \rangle \alpha \\
\langle x , y \rangle = \overline{\langle y, x \rangle}, \\
\langle x , x \rangle = 0 \Leftrightarrow x = 0.
\end{cases}
\end{equation}

We consider the norm $\| \cdots \|$ on $H$, defined by:
\begin{equation}\label{eq:2.12}
\|x\| = \langle x , x \rangle^{\frac{1}{2}}, \quad {}^\forall x \in H.
\end{equation}

Furthermore, if $(H, \| \cdot \|)$ is a complete normed space, then it will be called a \underline{right Hamilton space}.
\end{definition}

We will omit the proofs of the following Theorem~\ref{theorem:2.4} and the subsequent Corollaries, since we can mmick that of Theorem~\ref{theorem:2.1} and resort to elementary arguments

\begin{theorem}\label{theorem:2.4}
Every bounded linear operator on a separable infinite dimensional right Hamilton space has a nontrivial closed invariant subspace for it.
\end{theorem}

\begin{corollary}\label{corollary:2.5}
Let $T$ be a bounded linear operator on a separable infinite dimensional Hilbert space over $\mathbb{R}$ or $\mathbb{C}$ or right Hamilton space $H$. Let $X$ be a closed invariant subspace of $H$ for $T$ such that
\begin{equation}\label{eq:2.13}
2 \leq \textup{dim} (H \ominus X).
\end{equation}

Then $T$ has a nontrivial closed invariant subspace $Y$ of $H$ for $T$ such that
\begin{equation}\label{eq:2.14}
X \subsetneq Y \subsetneq H. \quad ([12])
\end{equation}
\end{corollary}

\begin{corollary}\label{corollary:2.6}
Let $H$ be a separable infinite dimensional Hilbert or right Hamilton space. Let $T$ be an injective closed linear operator in $H$ such that
\begin{equation}\label{eq:2.15}
\overline{\mathcal{D} (T)} = H,
\end{equation}
where $\mathcal{D}(T)$ denotes the domain of $T$, and
\begin{equation}\label{eq:2.16}
T ( \mathcal{D}(T)) \subset \mathcal{D}(T),
\end{equation}
admitting a unit cyclic vector $e_1 \in \mathcal{D}(T)$ in the sense that
\begin{equation}\label{eq:2.17}
H = [ T^{n-1} e_1: n \in \mathbb{N}].
\end{equation}

Further assume that
\begin{equation}\label{eq:2.18}
\sum_{n \in \mathbb{N}} \| T^{n-1} x \|^2 < \infty, \quad {}^\forall x \in \mathcal{D}(T).
\end{equation}

Then, by defining the new norm $\| \cdot \| $ on $\mathcal{D}(T)$ by:
\begin{equation}\label{eq:2.19}
\vertiii{x} = \big( \sum_{n \in \mathbb{N}} \| T^{n-1} x \|^2 \big)^{\frac{1}{2}}, \quad {}^\forall x \in \mathcal{D}(T),
\end{equation}
$T$ has a nontrivial $\vertiii{\cdot}$-closed invariant subspace of $(\mathcal{D}(T), \vertiii{\cdot})$ for $T$. ([8] Theorem~2.5)
\end{corollary}

\section*{Acknowledgements} 
I thank deeply Prof. Charles Akemann (UCSB), Prof. Sang Hoon Lee (CNU), Heon Lee (SNU), and my wife Soon Hee Kim.





\begin{center}
\line(1,0){350}
\end{center}
\end{document}